\documentclass[12pt]{amsart}

\usepackage{graphicx,colordvi,amssymb,colortbl,color,epsfig}
\usepackage{enumerate}
\usepackage[top=1in, bottom=1in, left=1in, right=1in]{geometry}
\usepackage{algorithm}
\usepackage{algpseudocode}

\newcommand{\dia}{$\diamond$}

\newcommand{\zpk}{\Z/\!\left(p^k\right)}

\newcommand{\Z}{\mathbb{Z}}
\newcommand{\N}{\mathbb{N}}
\newcommand{\cT}{\mathcal{T}}

\newcommand{\tf}{\tilde{f}}

\newcommand{\Q}{\mathbb{Q}}

\newcommand{\F}{\mathbb{F}}

\newcommand{\eps}{\varepsilon}
\newcommand{\ord}{\mathrm{ord}}

\newtheorem{dfn}{Definition}[section]
\newtheorem{lemma}[dfn]{Lemma}

\newtheorem{rem}[dfn]{Remark}
\newtheorem{ex}[dfn]{Example}
\newtheorem{thm}[dfn]{Theorem}

\newtheorem{algor}[dfn]{Algorithm}
\renewcommand{\qed}{$\blacksquare$}

\newcommand{\floor}[1]{\left\lfloor #1 \right\rfloor}
\newcommand{\ceil}[1]{\left\lceil #1 \right\rceil}

\newcommand{\rmv}[1]{}

\author{Leann Kopp} 
\email{LHK0002@auburn.edu}
\address{Mathematics Department, 
Auburn University, 221 Parker Hall, Auburn, Alabama \ 36849} 
\author{Natalie Randall} 
\email{natgrandall@gmail.com} 
\address{Department of Mathematics, 
Austin College, 900 N. Grand Ave.
Sherman, TX \ 75090} 
\author{J.\ Maurice Rojas}
\email{rojas@math.tamu.edu} 
\address{TAMU 3368, College Station, TX \ 77843-3368} 
\thanks{Partially supported by NSF grant CCF-1409020, and 
NSF REU grants DMS-1757872 and DMS-1460766. }  
\author{Yuyu Zhu} 
\email{zhuyuyu@math.tamu.edu}  
\address{TAMU 3368, College Station, TX \ 77843-3368 }

\title{\mbox{}\\ 
\vspace{-1in} 
Randomized Polynomial-Time Root Counting in Prime Power Rings}
\begin{document}

\begin{abstract} 
Suppose $k,p\!\in\!\N$ with $p$ prime and 
$f\!\in\!\Z[x]$ is a univariate polynomial with degree $d$ and 
all coefficients having absolute value less than $p^k$. 
We give a Las Vegas randomized 
algorithm that computes the number of roots of $f$ in $\zpk$ 
within time $\displaystyle{d^3(k\log p)^{2+o(1)}}$.  
(We in fact prove a more intricate complexity bound that is slightly 
better.) The best previous general algorithm had (deterministic)  
complexity exponential in $k$.  
We also present some experimental data evincing the potential practicality 
of our algorithm. 
\end{abstract} 

\maketitle

\pagenumbering{arabic} 

\section{Introduction}
Suppose $k,p\!\in\!\N$ with $p$ prime and $f\!\in\!\Z[x]$ is a univariate 
polynomial with degree $d\!\geq\!1$ and all coefficients having absolute 
value less than $p^k$. Let $N_{p,k}(f)$ denote the number of roots of $f$ 
in $\zpk$ (see, e.g., \cite{ragha,mcdo,bs,klivans,vzg,sircana} for further 
background on prime power rings). Computing $N_{p,k}(f)$ is a fundamental 
problem occuring 
in polynomial factoring \cite{lll,chistov,cantorqp,salagean,gnp}, coding theory 
\cite{lecerf}, and cryptography \cite{lauder}. The function $N_{p,k}(f)$ 
is also a basic ingredient in the study of Igusa zeta functions 
and algorithms over $\Q_p$ 
\cite{igusa,denef,cohenant,cantorqp,zg,denefver,lauderwan,dwan,chambert,airr}. 

In spite of the fundamental nature of computing $N_{p,k}(f)$, the 
fastest earlier general algorithms had complexity 
exponential in $k$: \cite{cgrw2} gave a deterministic algorithm taking time 
$\left(d\log(p) + 2^k\right)^{O(1)}$. 
While the $O$-constant was not stated in 
\cite{cgrw2}, the proof of the main theorem there indicates that the 
dependence on $k$ in their algorithm is linear in $e^k$. 
Note that counting the roots via brute-force takes time 
$dp^k(k\log p)^{1+o(1)}$, so the algorithm from \cite{cgrw2} is 
preferable, at least theoretically, for $p\!\geq\!3$. 
Here, we present a simpler, dramatically faster randomized algorithm 
(Algorithm \ref{algor:main} of the next section) that appears practical for 
all $p$. 
\begin{thm}
\label{thm:tree} 
Following the notation above, there is a Las Vegas randomized algorithm that 
computes $N_{p,k}(f)$ in time  
$\displaystyle{ 
kd^3(k\log p)^{1+o(1)}+\left(dk\log^2 p\right)^{1+o(1)} 
}$. In particular, the number of random bits needed is $O(dk\log(dk)\log p)$, 
and the space needed is $O(dk^2\log p)$ bits.
\end{thm} 

\noindent 
We prove Theorem \ref{thm:tree} in Section \ref{sec:proofs} below. 
In our context, {\em Las Vegas randomized} means that, with a 
fixed error probability (which we can take to be, say,
$\frac{1}{3}$), our algorithm under-estimates the number of roots.  
Our algorithm otherwise gives a correct root count, {\em and} always 
correctly announces whether the output count is correct or not. This 
type of randomization is standard in many number-theoretic algorithms, 
such as the fastest current algorithms for factoring polynomials over finite 
fields or primality checking (see, e.g., \cite{bs,kedlaya,cheng}).

At a high level, our algorithm here and the algorithm from 
\cite{cgrw2} are similar in that they reduce 
the main problem to a collection of computations, mostly in the finite field 
$\Z/(p)$, indexed by the nodes of a tree of size depending on 
$f$ and $k$. Also, both algorithms count by partitioning the roots in 
$\zpk$ into 
clusters having the same mod $p$ reduction. 
One subtlety to be aware of is that we compute the 
{\em number} of roots in $\zpk$, {\em without} listing all of them.  
Indeed, the number of roots in $\zpk$ can be as high as, say,  
$p^{k-\ceil{k/2}}$ (when $k\!\geq\!d\!=\!2$) or $p^{d/p}$ (when $d\!=\!kp$):  
simply consider the polynomials $x^2$ and $(x^p-x)^{d/p}$.   
So we can't attain a time or space bound sub-exponential in $k$ unless we do 
something more clever than naively store every root (see Remark \ref{rem:vzg} 
below).  

In finer detail, the algorithm from \cite{cgrw2} solves a ``small'' 
polynomial system at each node of a recursion tree (using a specially 
tailored Gr\"obner basis computation \cite{GVW}), while our algorithm performs 
a univariate factorization in $(\Z/(p))[x]$ at each node of a smaller 
recursion tree. Our use of fast factorization 
(as in \cite{kedlaya}) is why we avail to  
randomness, but this pays off: Gaining access to individual roots in 
$\Z/(p)$ (as suggested in \cite{cgrw2}) enables us to give a more 
streamlined algorithm. 

\begin{rem} 
\label{rem:vzg}  
von zur Gathen and Hartlieb presented in \cite{vzg} a randomized 
polynomial-time algorithm 
to compute all factorizations of certain $f\!\in\!(\Z/(p^k))[x]$. 
(Examples like $x^2\!=\!(x-p)(x+p)\!\in\!(\Z/(p^2))[x]$ show that 
unique factorization fails badly for $k\!\geq\!2$, and the number of  
possible factorizations can be exponential in $k$.) Their algorithm 
is particularly interesting since it uses a compact data structure to 
encode all the (possibly exponentially many) factorizations of $f$. 
Unfortunately, their algorithm has the restriction that $p^k$ not 
divide the discriminant of $f$. Their complexity bound, in our 
notation, is the sum of $d^7k\log(p)(k\log(p)+\log d)^2$ and a 
term involving the complexity of finding the mod $p^k$ reduction of a 
factorization over $\Z_p[x]$ (see Remarks 4.10--4.12 from \cite{vzg}). The 
complexity of just counting the number of possible factorizations (or just the
number of possible linear factors) of $f$  
from their data structure does not appear to be stated. \dia   
\end{rem} 

\noindent 
Creating an efficient classification of the roots of 
$f$ in $\zpk$ (and improving the data structure from \cite{vzg} 
by removing all restrictions on $f$), within time polynomial 
in $d+k\log p$, is a problem we hope to address in future work.  

For the reader interested in implementations, we have a preliminary 
{\tt Maple} implementation of Algorithm \ref{algor:main} freely downloadable 
from {\tt www.math.tamu.edu/\~{}rojas/count.map} . A few timings 
(all done on a Dell XPS13 Laptop with 8Gb RAM and a 256Gb ssd,  
running {\tt Maple 2015} within Ubuntu Linux 14.04) are listed below:

\noindent 
\mbox{}\hfill 
\begin{tabular}{cccc} 
$f(x)$  & $p^k$ & Brute-force$^1$  & Algorithm \ref{algor:main} \\ 
\hline 
Random degree $15$ & $2^{250}$ & $\approx\!2\times 10^{62}$ years & 0.077sec. \\ 
Random degree $75$ & $10009^{15}$ & $\approx\!5\times 10^{47}$ years & 0.116sec. \\ 
$(x-1234)^3(x-7193)^4(x-2030)^{12}$ & $123456791^1$ & 
 9min.\ 18sec. & 20.075sec \\ 
$(x-1234)^3(x-7193)^4(x-2030)^{12}$ & $123456791^{23}$ & 
 $\approx\!10^{173}$ years & 1min.\ 50.323sec.\\ \hline  
\end{tabular} \addtocounter{footnote}{1} 
\footnotetext{The timings
in years were based on extrapolating (without counting the
necessary expansion of laptop memory beyond 8Gb) from examples with much
smaller $k$ already taking over an hour.}
\hfill\mbox{}

\noindent 
Our {\tt Maple} implementations of brute-force and Algorithm \ref{algor:main} 
here are 5 lines long and 16 lines long, respectively. 
In particular, our random $f$ above were generated by 
taking uniformly random integer coefficients in $\{0,\ldots,p^k-1\}$ 
and then multiplying $5$ (or $25$) random cubic examples together: This 
results in longer timings for our code than directly picking a single random  
polynomial of high degree. The actual numbers of roots in the last two 
examples are respectively $3$ and\\  
\scalebox{.9}[1]{$83524650739763670783591272793501499347381420700990366689774050080031654011699848668752654$}\\
\scalebox{.9}[1]{$473531540039924209209663876325122031629580404523246324540823308088725469492593973$}. 

\subsection{A Recurrence from Partial Factorizations}  
Throughout this paper, we will use the integers $\{0,\ldots,p^k-1\}$ 
to represent elements of $\zpk$, unless otherwise specified. 
With this understanding, we will use the following notation: 
\begin{dfn} 
\label{dfn:basic} 
For any $f\!\in\!\Z[x]$ we let $\tf$ denote the mod $p$ reduction of $f$ 
and, for any root $\zeta_0\!\in\!\{0,\ldots,p-1\}$ of $\tf$,  
we call $\zeta_0$ {\em degenerate} if and only if $f'(\zeta_0)\!=\!0$ 
mod $p$. Letting $\ord_p: \Z \longrightarrow \N\cup\{0\}$ denote the 
usual $p$-adic valuation with $\ord_p(p)\!=\!1$, we then define 
$s(f,\eps)\!:=\!\min\limits_{j\geq 0}\!\left\{j+\ord_p 
\frac{f^{(j)}(\eps)}{j!}\right\}$ for any $\eps\!\in\!\{0,\ldots,p-1\}$. 
Finally, fixing $k\!\in\!\N$, let us inductively define a set 
$T_{p,k}(f)$ of pairs $(f_{i,\zeta},k_{i,\zeta})\!\in\!\Z[x]\times \N$ as 
follows: We set $(f_{0,0},k_{0,0})\!:=\!(f,k)$. Then, for 
any $i\!\geq\!1$ with $(f_{i-1,\mu},k_{i-1,\mu})\!\in\!T_{p,k}(f)$ and 
any {\em degenerate} root $\zeta_{i-1}\!\in\!\{0,\ldots,p-1\}$ of 
$\tilde{f}_{i-1,\mu}$ with 
$s_{i-1}\!:=\!s(f_{i-1,\mu},\zeta_{i-1})\!\in\!\{2,\ldots,k_{i-1,\mu}-1\}$, we
define $\zeta\!:=\!\mu+p^{i-1}\zeta_{i-1}$, $k_{i,\zeta}\!:=\!k_{i-1,\mu}
-s_{i-1}$ and $f_{i,\zeta}(x)\!:=\! 
\left[\frac{1}{p^{s_{i-1}}}f_{i-1,\mu}(\zeta_{i-1}+px)\right] \ \mathrm{mod} 
\ p^{k_{i,\zeta}}$. \dia 
\end{dfn} 

\noindent 
The ``perturbations'' $f_{i,\zeta}$ of $f$ will help us keep 
track of how the roots of $f$ in $\zpk$ cluster (in a $p$-adic 
metric sense) about the roots of $\tf$. 
Since $\frac{f^{(j)}(\eps)}{j!}$ is merely the coefficient of 
$x^j$ in the Taylor expansion of $f(x+\eps)$ about $x\!=\!0$,  
it is clear that $\frac{f^{(j)}(\eps)}{j!}$ 
is always an integer (under the assumptions above) provided 
$\eps\!\in\!\Z$. 

We will see in the next section how $T_{p,k}(f)$ can be identified with a 
finite rooted directed tree. In particular, it is easy to see that the set 
$T_{p,k}(f)$ is always finite since, by construction, only 
$f_{i,\zeta}$ with $i\!\leq\!\floor{(k-1)/2}$ and $\zeta\!\in\!\Z/(p)$ 
are possible (see also Lemma \ref{lemma:branch} of Section \ref{sec:proofs} 
below). 
\begin{ex}
\label{ex:tri}
Let us take $p\!=\!3$, $k\!=\!7$, and $f(x)\!:=\!x^{10}-10x+738$. 
A simple calculation then shows that 
$\tf_{0,0}(x)\!=\!x(x-1)^9$, which has roots $\{0,1\}$ in $\Z/(3)$. 
The root $0$ is non-degenerate so the only possible $f_{1,\zeta}$ 
would be an $f_{1,1}\!=\!f_{1,0+1}$. 

In particular, $s(f_{0,0},1)\!=\!4$ and thus $k_{1,1}\!=\!3$ and 
$f_{1,1}(x)\!=\!21x^4+13x^3+5x^2+9$ mod $3^3$.  
Since $\tf_{1,1}(x)\!=\!x^2(x-1)$ and $1$ is a non-degenerate
root of $\tf_{1,1}$, we see that the only possible 
$f_{2,\zeta}$ would be an $f_{2,1}\!=\!f_{2,1+0}$.  

Since $s(f_{1,1},0)\!=\!2$ we then obtain $k_{2,1}\!=\!1$, and 
$f_{2,1}(x)\!=\!2(x-1)(x-2)$ mod $3$, which has only non-degenerate roots. 
So by Definition \ref{dfn:basic} there can be no $f_{3,\zeta}$ and thus 
our collection of pairs $T_{p,k}(f)$ consists of just $3$ pairs. \dia
\end{ex}

Using base-$p$ expansion, there is an obvious bijection between the 
ring $\Z_p$ of $p$-adic integers and the set of root-based paths in 
an infinite $p$-ary tree $\mathbb{T}_p$. It is then natural to use the 
leafs of a finite subtree of $\mathbb{T}_p$ to 
store the roots of $f$ in $\zpk$. This type of tree structure was studied 
earlier by Schmidt and Stewart in \cite{schmidt,schmidtstewart}, from the 
point of view of classification and (in our notation) upper bounds on 
$N_{p,k}(f)$. However, it will be more algorithmically efficient to instead 
endow our set $T_{p,k}(f)$ with a tree structure. The following 
fundamental lemma relates $N_{p,k}(f)$ to a recursion tree structure on 
$T_{p,k}(f)$. 
\begin{lemma} 
\label{lemma:basic} 
Following the notation above, let $n_p(f_{0,0})$ denote the number of 
non-degenerate roots of $\tf_{0,0}$ in $\Z/(p)$. 
Then, provided $k\!\geq\!2$ and $f_{0,0}$ is not identically $0$ in 
$(\Z/(p))[x]$, we have  
\[N_{p,k}(f_{0,0}) \ = \ n_p(f_{0,0}) 
+ \left(\sum\limits_{\substack{\zeta_0\in\Z/(p)\\
                        s(f_{0,0},\zeta_0)\geq k}} p^{k-1}\right) 
+ \!\!\!\!\!\!\! 
\sum\limits_{\substack{\zeta_0\in\Z/(p)\\
                        s(f_{0,0},\zeta_0)\in\{2,\ldots,k-1\}}} 
\!\!\!\!\!\!\! 
p^{s(f_{0,0},\zeta_0)-1}N_{p,k-s(f_{0,0},\zeta_0)}(f_{1,\zeta_0})
.\]  
\end{lemma}  

\noindent 
We prove Lemma \ref{lemma:basic} in the next section, where 
it will immediately follow that Lemma \ref{lemma:basic} applies recursively, 
i.e., our root counting formula still holds if one replaces  
$(f_{0,0},k,f_{1,\zeta_0},\zeta_0)$ with 
$(f_{i-1,\mu},k_{i-1,\mu},f_{i,\mu+p^{i-1}\zeta_{i-1}},\zeta_{i-1})$. 
There we also show how Lemma \ref{lemma:basic} leads to our 
recursive algorithm (Algorithm \ref{algor:main}) 
for computing $N_{p,k}(f)$. In essence, the third sum term 
above is what creates children for a node corresponding 
to $f_{i-1,\mu}$. 

Note that by construction, 
$s(f,\zeta_0)\!\geq\!2$ implies that $\zeta_0$ is a degenerate root of 
$\tf$. So the last two sum terms in the formula (from Lemma \ref{lemma:basic} 
above) range over certain degenerate roots of 
$\tf$. Note also that $N_{p,k}(f)$ depends only on the 
residue class of $f$ mod $p^k$, so we will often abuse the notations 
$N_{p,k}(f)$ and $s(f,\zeta_0)$ by allowing 
$f\!\in\!\left(\zpk\right)[x]$ as well. The following 
example illustrates how $N_{p,k}(f)$ can be computed recursively. 

\medskip 
\noindent 
{\bf Example 1.6.} \addtocounter{dfn}{1} 

\vspace{-.3cm} 
\noindent 
\begin{minipage}[t]{.7 \textwidth}
\vspace{0pt}
{\em Revisiting Example \ref{ex:tri}, let us 
count the roots in $\Z/\!\left(3^7\right)$ of\linebreak 
$f(x)\!:=\!x^{10}-10x+738$. 
Lemma \ref{lemma:basic} and our earlier computation of $T_{p,k}(f)$ 
then tell us that $N_{3,7}(f)\!=\!1+3^3 N_{3,3}(f_{1,1})$ and 
$N_{3,3}(f_{1,1})\!=\!1+3^1 N_{3,1}(f_{2,1})$ where $f_{2,1}(x)\!=\!
2(x-1)(x-2)$. So we obtain 
$N_{3,7}(f)\!=\!1+3^3(1+3^1\cdot 2)\!=\!190$. (Our 
{\tt Maple} implementation confirmed this count in under $4$ milliseconds.) 
We illustrate the corresponding tree structure (defined in Section 
\ref{sec:proofs} below) on the right. Note that the powers of $3$ in the 
expression $1+3^3(1+3^1\cdot 2)$ appear as edge labels in our tree, but 
the contribution of {\em non}-degenerate roots to our count 
is not notated on our tree. \dia}
\end{minipage}   
\begin{minipage}[t]{.3 \textwidth}
\vspace{0pt}
\scalebox{.8}[.8]{
\vbox{
\begin{picture}(200,0)(20,305)
\put(90,116){\epsfig{file=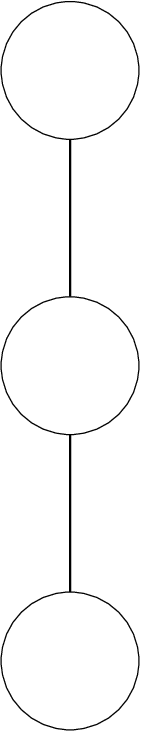,height=3.1in,clip=}}
\put(92,316){$(f_{0,0},7)$}     
  \put(116,273){$3^3$}  
\put(92,225){$(f_{1,1},3)$} 
  \put(116,184){$3^1$}
\put(92,134){$(f_{2,1},1)$} 
\end{picture}}}
\end{minipage}

\vspace{.3cm} 
\noindent 
While the tree from our example above has just $3$ nodes, 
the earlier tree structure from \cite{schmidt,schmidtstewart}  
would have resulted in over $190$ nodes. We will now fully detail 
how to efficiently reduce root counting over $\zpk$ to 
computing $p$-adic valuations and factoring in $(\Z/(p))[x]$. 

\section{Algebraic Preliminaries and Our Algorithm} 
\label{sec:prelim} 
Let us first recall the following version of Hensel's Lemma: 
\begin{lemma}(See, e.g., \cite[Thm.\ 2.3, Pg.\ 87, Sec.\ 2.6]{nzm}.)  
\label{lemma:hensel} 
Suppose $k\!\in\!\N$, $f\!\in\!\Z[x]$ is not identically $0$ in  
$(\Z/(p))[x]$, and $\zeta_0\!\in\!\Z/(p)$ is a non-degenerate root of $\tf$. 
Then there is a unique $\zeta\!\in\!\zpk$ with $\zeta\!=\!\zeta_0$ 
$\mathrm{mod} \ p$ and $f(\zeta)\!=\!0 \ \mathrm{mod}  \ p^k$. \qed 
\end{lemma} 

The following lemma enables us to understand the lifts of degenerate roots 
of $\tf$. 
\begin{lemma} 
\label{lemma:mult} 
Following the notation of Lemma \ref{lemma:hensel}, 
suppose instead that $\zeta_0\!\in\!\Z/(p)$ is a root of 
$\tf$ of (finite) multiplicity $m\!\geq\!2$. Suppose also that $k\!\geq\!2$ 
and that there is a $\zeta\!\in\!\zpk$ with $\zeta\!=\!\zeta_0$ mod $p$ and 
$f(\zeta)\!=\!0$ mod $p^k$. Then $s(f,\zeta_0)\!\in\!\{2,\ldots,m\}$. 
\end{lemma} 

\noindent 
{\bf Proof of Lemma \ref{lemma:mult}:} We may assume, by base-$p$ 
expansion that $\zeta\!=\!
\zeta_0+p\zeta_1+\cdots+p^{k-1}\zeta_{k-1}$ for some $\zeta_1,
\ldots,\zeta_{k-1}\!\in\!\{0,\ldots,p-1\}$. Note that 
$f'(\zeta_0)\!=\!0$ mod $p$ since $\zeta_0$ is a degenerate root. 
Note also that $j+\ord_p\frac{f^{(j)}(\zeta_0)}{j!}\!\geq\!2$ for all 
$j\!\geq\!2$. Letting 
$\sigma\!:=\!\zeta_1+p\zeta_2+\cdots+p^{k-2}\zeta_{k-1}$ 
we then see by Taylor expansion that 
$f(\zeta)\!=\!f(\zeta_0)+f'(\zeta_0)p\sigma+\cdots
+\frac{f^{(k-1)}}{(k-1)!}p^{k-1}\sigma^{k-1}$ mod $p^k$. So $f(\zeta)\!=\!0$ 
mod $p^k$ implies that $f(\zeta_0)\!=\!0$ mod $p^2$ and thus 
$s(f,\zeta_0)\!\geq\!2$.

To conclude, by Taylor expansion about $\zeta_0$, our multiplicity assumption 
implies that\linebreak 
\scalebox{1}[1]{$\frac{f^{(m)}(\zeta_0)}{m!}$ is an integer not divisible by 
$p$. So $m+\ord_p\frac{f^{(m)}(\zeta_0)}{m!}\!=\!m$ and thus 
$s(f,\zeta_0)\!\leq\!m$. \qed}   

\medskip 
We are now ready to state our main algorithm. 

\noindent 
\mbox{}\hspace{.7cm}\scalebox{.87}[.87]{\fbox{\mbox{}\hspace{.3cm}\vbox{
\begin{algor}[{\tt RandomizedPrimePowerRootCounting}$(f,p,k)$] 
\label{algor:main}
\mbox{}\\ 
{\bf Input.} 
$(f,p,k)\!\in\!\Z[x]\times \N\times \N$ with 
$p$ prime and $f(x)\!=\!c_0+\cdots+c_dx^d$.\\ 
{\bf Output.} An integer $M\!\leq\!N_{p,k}(f)$ 
that, with probability at least $\frac{2}{3}$, is exactly $N_{p,k}(f)$.\\ 
{\bf Description.} \\ 
1: {\bf Let} $v\!:=\!\min\limits_{i\in \{0,\ldots,d\}}\ord_p c_i$ and 
    $f_{0,0}\!:=\!f$.\\  
2: {\bf If} $v\!\geq\!k$ \\ 
3: \mbox{}\hspace{.7cm}{\bf Let} $M\!:=\!p^k$. {\bf Return}.\\   
4: {\bf Elseif} $v\!\in\!\{1,\ldots,k-1\}$\\ 
5: \mbox{}\hspace{.7cm}{\bf Let} $M\!:=\!p^v\text{{\tt 
RandomizedPrimePowerRootCounting}$\left(\frac{f_{0,0}(x)}{p^v},p,k-v\right)$}$. 
{\bf Return}.\\ 
6: {\bf End(If)}. \\ 
7: {\bf Let} $M\!:=\!\deg \gcd(\tf_{0,0},x^p-x)$. \\  
8: {\bf For} $\zeta_0\!\in\!\Z/(p)$ a degenerate root of $\tf_{0,0}$ 
     {\bf do}$^2$\\ 
9: \mbox{}\hspace{.7cm}{\bf Let} $s\!:=\!s(f_{0,0},\zeta_0)$.\\   
\mbox{}\hspace{-.2cm}10: \mbox{}\hspace{.7cm}{\bf If} $s\!\geq\!k$\\ 
\mbox{}\hspace{-.2cm}11: \mbox{}\hspace{1.4cm}{\bf Let} $M\!:=\!M+p^{k-1}$.\\ 
\mbox{}\hspace{-.2cm}12: \mbox{}\hspace{.7cm}{\bf Elseif} $s\!\in\{2,\ldots,
k-1\}$\\ 
\mbox{}\hspace{-.2cm}13: \mbox{}\hspace{1.4cm}{\bf Let} 
$M\!:=\!M+p^{s-1}\text{{\tt 
RandomizedPrimePowerRootCounting}$\left(f_{1,\zeta_0},p,k-s\right)$}$.\\ 
\mbox{}\hspace{-.2cm}14: \mbox{}\hspace{.7cm}{\bf End(If)}.\\   
\mbox{}\hspace{-.2cm}15: {\bf End(For)}.\\ 
\mbox{}\hspace{-.2cm}16: {\bf If} the preceding {\bf For} loop did {\em not}  
access all the degenerate roots of $\tf_{0,0}$ \\ 
\mbox{}\hspace{-.2cm}17: \mbox{}\hspace{.7cm} Print {\tt ``Sorry, 
your Las Vegas factoring method failed.}\\ 
\mbox{}\hspace{2.9cm}{\tt You have an under-count so you should try 
re-running.''}\\ 
\mbox{}\hspace{-.2cm}18: {\bf End(If)}. \\ 
\mbox{}\hspace{-.2cm}19: Print {\tt ``If you've seen no under-count messages 
then your count is correct!''}\\  
\mbox{}\hspace{-.2cm}20: {\bf Return}. 
\end{algor}}
}
}
\addtocounter{footnote}{1} 
\footnotetext{Here we use the fastest available Las Vegas factoring 
algorithm over $(\Z/(p))[x]$ (currently \cite{kedlaya}) to 
isolate the degenerate roots of $\tf$. Such factoring algorithms 
enable us to correctly announce failure to find all the degenerate 
roots, should this occur. We describe in the next section how to 
efficiently control the error probability. } \\ 

Before proving the correctness of Algorithm \ref{algor:main}, it will 
be important to prove our earlier key lemma. 

\medskip 
\noindent 
{\bf Proof of Lemma \ref{lemma:basic}:} 
Proving our formula clearly reduces to determining how many 
lifts each possible root $\zeta_0\!\in\!\Z/(p)$ of $\tf_{0,0}$ has 
to a root of $f_{0,0}$ in $\zpk$. Toward this end, note that 
Lemma \ref{lemma:hensel} implies that each non-degenerate $\zeta_0$ 
lifts to a unique root of $f_{0,0}$ in $\zpk$. In particular, this 
accounts for the summand $n_p(f_{0,0})$ in our formula. So now we merely 
need to count the lifts of the degenerate roots. 

Assume $\zeta_0\!\in\!\Z/(p)$ is a degenerate root of $\tf_{0,0}$, 
write $\zeta\!=\!\zeta_0+p\zeta_1+\cdots+p^{k-1}\zeta_{k-1}\!\in\!\zpk$ 
via base-$p$ expansion as before, set $\sigma\!:=\!
\zeta_1+p\zeta_2+\cdots+p^{k-2}\zeta_{k-1}$, and let $s\!:=\!s(f_{0,0},
\zeta_0)$. Clearly then, $f_{0,0}(\zeta)\!=\!p^{s}f_{1,\zeta_0}(\sigma)$ mod 
$p^k$ and, by construction, $f_{1,\zeta_0}\!\in\!\Z[x]$ and is not identically 
$0$ in $(\Z/(p))[x]$. 

If $s\!\geq\!k$ then $f_{0,0}(\zeta)\!=\!0$ mod $p^k$ independent of $\sigma$. 
So there are exactly $p^{k-1}$ values of $\zeta\!\in\!\zpk$ with 
$\zeta\!=\!\zeta_0$ mod $p$. This accounts for the second summand in our 
formula. 

If $s\!\leq\!k-1$ then $\zeta$ is a root of $f_{0,0}$ 
with $\zeta\!=\!\zeta_0$ mod $p$ if and only if $f_{1,\zeta_0}(\sigma)\!=\!0$ 
mod $p^{k-s}$. Also, $s\!\geq\!2$ (thanks to Lemma \ref{lemma:mult}) 
because $\zeta_0$ is a degenerate root. 
Since the base-$p$ digits $\zeta_{k-s+1},\ldots,\zeta_{k-1}$
do not appear in the last equality, the number of possible lifts
$\zeta$ of $\zeta_0$ is thus\linebreak  

\vbox{
\noindent 
exactly $p^{s-1}$ times the number of 
roots $\zeta_1+p\zeta_2+\cdots+p^{k-s-1}\zeta_{k-s}\!\in\!\Z/\!\left(p^{k-s}
\right)$ of $f_{1,\zeta_0}$. So this 
accounts for the third summand in our formula and we are done. \qed 

\medskip 
We are at last ready to prove the correctness of Algorithm \ref{algor:main}. 

\medskip 
\noindent 
{\bf Proof of Correctness of Algorithm \ref{algor:main}:} 
Assume temporarily that Algorithm \ref{algor:main} is correct 
when $f_{0,0}$ is {\em not} identically $0$ in $(\Z/(p))[x]$. 
Since (for any integers $a,x,y$ with $a\!\leq\!k$)  
$p^a x\!=\!p^a y$ mod $p^k \Longleftrightarrow 
x\!=\!y$ mod $p^{k-a}$, Steps 1--6 of our algorithm then clearly 
correctly dispose of the case where $f$ is identically $0$ in 
$(\Z/(p))[x]$. So let us now prove correctness when $f$ is {\em not} 
identically $0$ in $(\Z/(p))[x]$. Applying Lemma \ref{lemma:basic}, we then 
see that 
it is enough to prove that the value of $M$ is the value of our formula for 
$N_{p,k}(f)$ when the {\bf For} loops of Algorithm \ref{algor:main} runs 
correctly. 

Step 7 ensures that the value of $M$ is initialized as $n_p(f)$. 
Steps 8--15 (once the {\bf For} loop is completed) then simply 
add the second and third summands of our formula to 
$M$ thus ensuring that $M\!=\!N_{p,k}(f)$, {\em provided} the  
{\bf For} loop has run correctly, along with all the {\bf For} loops 
in the recursive calls to {\tt RandomizedPrimePowerRootCounting}. 
Should any of these {\bf For} loops run incorrectly, Steps 16--20 ensure 
that our algorithm correctly announces an under-count.$^3$ So we are done. \qed 

\section{Our Complexity Bound: Proving Theorem \ref{thm:tree}} 
\label{sec:proofs} 
Let us now introduce a tree structure on $T_{p,k}(f)$ that will enable our 
complexity analysis. 
\begin{dfn} 
\label{dfn:tree} 
Let us identify the elements of $T_{p,k}(f)$ with nodes of a labelled rooted 
directed tree $\cT_{p,k}(f)$ defined inductively as follows: 
\begin{enumerate}
\item{We set $f_{0,0}\!:=\!f$, $k_{0,0}\!:=\!k$, and let 
$(f_{0,0},k_{0,0})$ be the label of the root node of 
$\cT_{p,k}(f)$. }   
\item{The non-root nodes of $\cT_{p,k}(f)$ are uniquely labelled by each 
$(f_{i,\zeta},k_{i,\zeta})\!\in\!T_{p,k}(f)$ with $i\!\geq\!1$.} 
\item{There is an edge from node $(f_{i',\zeta'},k_{i',\zeta'})$ to 
node $(f_{i,\zeta},k_{i,\zeta})$ if and only if 
$i'\!=\!i-1$ and there is a degenerate root $\zeta_{i-1}\!\in\!\Z/(p)$ 
of $\tf_{i',\zeta'}$ with $s(f_{i',\zeta'},\zeta_{i-1})\!\in\!\{2,
\ldots,k_{i',\zeta'}-1\}$ and $\zeta\!=\!\zeta'+p^{i-1}\zeta_{i-1}\!\in\!
\Z/(p^i)$.}  
\item{The label of a directed edge from node $(f_{i',\zeta'},k_{i',\zeta'})$ 
to node $(f_{i,\zeta},k_{i,\zeta})$ is $p^{s\left(f_{i',\zeta'},
(\zeta-\zeta')/p^{i'}\right)-1}$.} 
\end{enumerate}  
In particular, the edges are labelled by powers of 
$p$ in $\{p^1,\ldots,p^{k-2}\}$, and the labels of 
the nodes lie in $\Z[x]\times \N$. \dia 
\end{dfn}   

\noindent 
{\bf Example 3.2.} \addtocounter{dfn}{1} 

\vspace{-.3cm} 
\noindent 
\begin{minipage}[t]{.7 \textwidth} 
\vspace{0pt} 
{\em Letting 
$g(x)\!:=\!x^5-8x^4+25x^3-38x^2+28x-8$, 
the tree $\cT_{17,100}(g)$ is drawn to the right:  
Note that $\cT_{17,100}(g)$ has depth $\floor{(100-1)/2}\!=\!49$ 
and exactly $1+\floor{(100-1)/2}+\floor{(100-1)/3}\!=\!83$ nodes. 
To count the roots of $g$ in $\Z/\!\left(17^{100}\right)$ one can
then easily calculate that
$g_{0,0}(x)\!=\!(x-1)^2(x-2)^3$,
$g_{1,1}(x)\!=\!x^2(4913x^3-867x^2+51x-1)$ and 
$g_{1,2}(x)\!=\!x^3(289x^2+34x+1)$. 
The last 
two polynomials have {\em no} nonzero roots mod $17$. 

A bit more computation then yields
$\tilde{g}_{i,1}(x)\!=\!-x^2$ for all $i\!\in\!\{1,\ldots,49\}$ and
$\tilde{g}_{j,2}(x)\!=\!x^3$ for all $j\!\in\!\{1,\ldots,33\}$.  
Also,\linebreak 
\scalebox{.98}[1]{$N_{17,2}(g_{49,1})\!=\!17$ by Lemma \ref{lemma:basic} and
$N_{17,1}(g_{33,2})\!=\!1$ trivially. So by}
Lemma \ref{lemma:basic}
once more, $g_{0,0}$ has exactly $17\cdot 17^{49}\cdot 1 +17^2\cdot (17^2)^{32}
\cdot 1\!=$\linebreak 
\scalebox{.95}[1]{$17^{50}+17^{66}$ roots in $\Z/\!\left(17^{100}\right)$. 
Expanded in base-$10$, this count is}  
} 
\end{minipage}  
\begin{minipage}[t]{.3 \textwidth} 
\vspace{0pt} 
\scalebox{.58}[.58]{
\vbox{
\begin{picture}(200,0)(160,290) 
\put(220,22){\epsfig{file=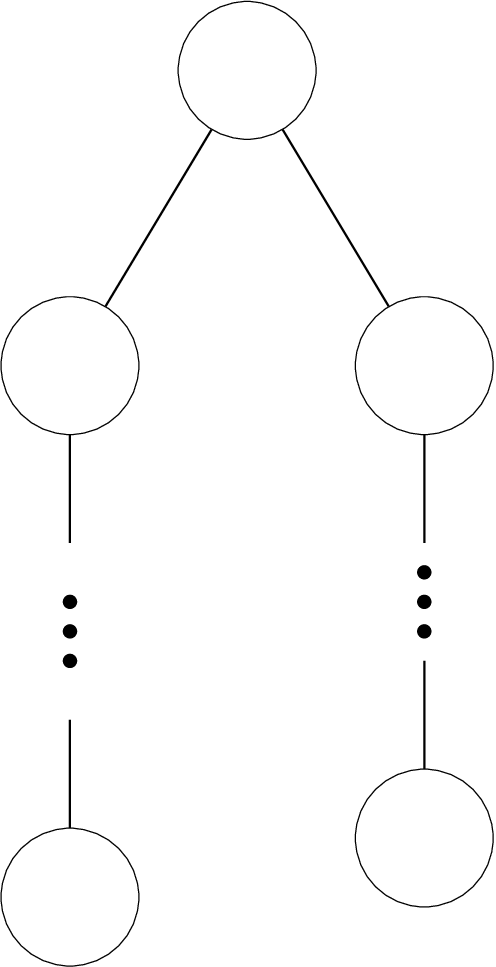,height=4.4in,clip=}} 
     \put(277,316){\scalebox{.9}[.9]{$(g_{0,0},100)$}} 
     \put(258,273){$17^1$}  \put(332,271){$17^2$}
     \put(224,219){\scalebox{.9}[.9]{$(g_{1,1},98)$}} 
     \put(340,219){\scalebox{.9}[.9]{$(g_{1,2},97)$}}
     \put(226,178){$17^1$} \put(361,179){$17^2$}
     \put(226,78){$17^1$}               \put(361,100){$17^2$}
     \put(224,42){\scalebox{.9}[.9]{$(g_{49,1},2)$}}     
     \put(341,63){\scalebox{.9}[.9]{$(g_{33,2},1)$}}
\end{picture}}}  
\end{minipage} 
\scalebox{.94}[1]{$1620424537653706124196923258781575759359875675913436470380245486276378993995166018$. \dia}   
}\addtocounter{footnote}{1} 
\footnotetext{Note that
checking for an under-count can be reduced to an irreducibility 
check in $\F_p[x]$, which can be done in deterministic polynomial-time: 
See, e.g., 
\cite[Cor.\ 14.35, Algor.\ 14.36, \& Thm.\ 14.37, pp.\ 406--408]{vzg}.} 

\begin{rem} 
Our trees $\cT_{p,k}(\cdot)$ thus encode algebraic expressions  
for our desired root counts $N_{p,k}(\cdot)$. In particular, 
the children of a node labelled $(f_i,k_i)$ 
yield terms (corresponding to the child nodes) that one sums to get 
the root count $N_{p,k_i}(f_i)$, and the edge labels yield weights  
multiplying the corresponding terms. (The contribution from 
non-degenerate roots is not visible from the tree but does 
influence each $N_{p,k_i}(f_i)$, as detailed by 
Lemma \ref{lemma:basic}.) \dia  
\end{rem} 
\begin{ex} 
Suppose we set $p\!=\!31$, $k\!=\!7$, and we define $h(x)$ to be\\  
\scalebox{.9}[1]{$x^{12}-60x^{11}-4420x^{10}+275040x^9+8287728x^8
-502626240x^7-8802489280 
x^6-10069291727x^5$}\\ 
\mbox{}\hspace{.5cm}\scalebox{.9}[1]{$-6168330858x^4-10982634616x^3
+6650045702x^2-4862117081x-6450915579$}.\\ 
Then the tree $\cT_{31,7}(h)$ has the following structure: \\  
\begin{picture}(200,150)(0,0)
\put(50,0){\epsfig{file=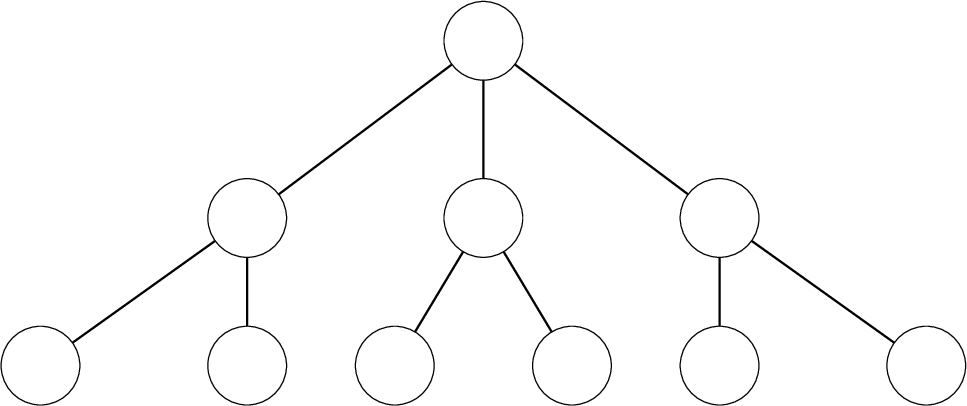,height=2in,clip=}}
\put(210,128){\scalebox{.6}[.6]{$(h_{0,0},7)$}} 
\put(173,103){\scalebox{.6}[.6]{$31^3$}} 
\put(210,95){\scalebox{.6}[.6]{$31^3$}} 
\put(263,103){\scalebox{.6}[.6]{$31^3$}} 
\put(126,65){\scalebox{.6}[.6]{$(h_{1,1},3)$}} 
\put(209,65){\scalebox{.6}[.6]{$(h_{1,15},3)$}} 
\put(293,65){\scalebox{.6}[.6]{$(h_{1,30},3)$}} 
\put(90,40){\scalebox{.6}[.6]{$31^1$}} 
\put(141,40){\scalebox{.6}[.6]{$31^1$}} 
\put(196,40){\scalebox{.6}[.6]{$31^1$}} 
\put(241,40){\scalebox{.6}[.6]{$31^1$}} 
\put(295,40){\scalebox{.6}[.6]{$31^1$}} 
\put(347,40){\scalebox{.6}[.6]{$31^1$}} 
\put(52,12){\scalebox{.6}[.6]{$(h_{2,32},1)$}} 
\put(124,12){\scalebox{.57}[.6]{$(h_{2,931},1)$}} 
\put(178,12){\scalebox{.6}[.6]{$(h_{2,46},1)$}} 
\put(240,12){\scalebox{.57}[.6]{$(h_{2,945},1)$}} 
\put(293,12){\scalebox{.6}[.6]{$(h_{2,61},1)$}} 
\put(366,12){\scalebox{.57}[.6]{$(h_{2,960},1)$}} 
\end{picture}

\noindent 
In particular, the polynomials corresponding to the depth $1$ nodes are 
exactly\\ 
\mbox{}\hspace{.5cm}$h_{1,1}\!=\!
9610x^6+13640x^5+25563x^4+2511x^3+9417x^2+13640x+14992$,\\  
\mbox{}\hspace{.5cm}$h_{1,15}\!=\!22103x^6+1674x^5+11825x^4+26443x^3+29205x^2
+1674x+26240$, and\\  
\mbox{}\hspace{.5cm}$h_{1,30}\!=\! 
24986x^6+28520x^5+22228x^4+2542x^3+29541x^2+28520x+12618$.\\  
Also, the polynomials corresponding to the depth $2$ nodes are exactly 
$h_{2,1+1\cdot 31}\!=\!h_{2,1+30\cdot 31}\!=\!14x^2$ and 
$h_{2,15+1\cdot 31}\!=\!h_{2,15+30\cdot 31}\!=\!h_{2,30+1\cdot 31}\!=\!
h_{2,30+30\cdot 31}\!=\!4x^2$. So Lemma \ref{lemma:basic} tells us that $h$ 
has exactly $6\cdot 31^4\cdot 1\!=\!5541126$ roots in $\Z/(31^7)$.   
\end{ex}

The following lemma will be central in our complexity analysis. 
\begin{lemma}
\label{lemma:ref}
Suppose $k,p\!\in\!\N$ with $p$ prime, $f\!\in\!\Z[x]$ has degree $d$, 
and $(f_{i-1,\mu},k_{i-1,\mu})$ is any node of 
$\cT_{p,k}(f)$. Then $\sum \deg \tilde{f}_{i,\zeta} 
\leq \deg \tilde{f}_{i-1,\mu}$, where the sum ranges over all 
child nodes $(f_{i,\zeta},k_{i,\zeta})$ of $(f_{i-1,\mu},k_{i-1,\mu})$. \qed 
\end{lemma}

\noindent 
Lemma \ref{lemma:ref} follows immediately from the last sentence of 
Assertion (3) of the more refined lemma below: 
\begin{lemma} 
\label{lemma:branch} 
Following the notation of Lemma \ref{lemma:ref}, we have that:  
\begin{enumerate} 
\item{The depth of $\cT_{p,k}(f)$ is at most $\floor{(k-1)/2}$.} 
\item{The degree of the root node of $\cT_{p,k}(f)$ is at most $\floor{d/2}$.} 
\item{The degree of any {\em non}-root node of $\cT_{p,k}(f)$ labelled 
$(f_{i,\zeta},k_{i,\zeta})$, with parent $(f_{i-1,\mu},k_{i-1,\mu})$ 
and \ $\zeta_{i-1}:=(\zeta-\mu)/p^{i-1}$, \ 
is \ at \ most \ $\floor{s(f_{i-1,\mu},\zeta_{i-1})/2}$. \ In \ 
particular,\\  
$\deg \tf_{i,\zeta}\!\leq\!s(f_{i-1,\mu},\zeta_{i-1})
\!\leq\!k_{i-1,\mu}-1\!\leq\!k-1$ and   
$\!\!\!\!\!
\!\!\!\!\!
\sum\limits_{\substack{(f_{i,\zeta},k_{i,\zeta}) \text{ a child}\\ 
\text{of } (f_{i-1,\mu},k_{i-1,\mu})}}  
\!\!\!\!\!
\!\!\!\!\!
s(f_{i-1,\mu},\zeta_{i-1})\leq \deg \tf_{i-1,\mu}$.  }
\item{$\cT_{p,k}(f)$ has at most $\floor{\frac{d}{2}}$ nodes at 
depth $i\!\geq\!1$, and thus a total of no more than 
$1+\floor{\frac{d}{2}}\floor{\frac{k-1}{2}}$ nodes.}  
\end{enumerate} 
\end{lemma} 

\vbox{
\noindent 
{\bf Proof of Lemma \ref{lemma:branch}:} 

\medskip 
\noindent 
{\bf Assertion (1):} By Definitions 
\ref{dfn:basic} and \ref{dfn:tree}, the labels $(f_{i,\zeta},k_{i,\zeta})$ 
satisfy\linebreak 
$2\!\leq\!k_{i-1,\mu}-k_{i,\zeta}\!\leq\!k_{i-1,\mu}-1$  
for any child $(f_{i,\zeta},k_{i,\zeta})$ of $(f_{i-1,\mu},k_{i-1,\mu})$, 
and $1\!\leq\!k_{i,\zeta}\!\leq\!k-2$ for all $i\!\geq\!1$. 
So considering any root to leaf path in $\cT_{p,k}(f)$,   
it is clear that the depth of $\cT_{p,k}(f)$ can be no greater than 
$1+\floor{(k-2-1)/2}\!=\!\floor{(k-1)/2}$. \qed } 

\medskip 
\noindent 
{\bf Assertion (2):} 
Since $\tf_{0,0}\!=\!\tf$ has degree $\leq\!d$, and the multiplicity of any 
degenerate root of $\tf_{0,0}$ is at least $2$, we see that $\tf_{0,0}$ has 
no more than $\floor{d/2}$ degenerate roots in $\Z/(p)$. Every edge 
emanating from the root node of $\cT_{p,k}(f)$ corresponds to a unique 
degenerate root of $\tf_{0,0}$ (and not every degenerate root of $\tf$ 
need yield a valid edge emanating from the root of $\cT_{p,k}(f)$), so we are 
done. \qed 

\medskip 
\noindent 
{\bf Assertion (3):} 
The degree bound for non-root nodes follows similarly to the 
degree bound for the root node:  
Letting $s\!:=\!s(f_{i-1,\mu},\zeta_{i-1})$, it 
suffices to prove that $\deg \tilde{f}_{i,\zeta}\!\leq\!s$ for all 
$i\!\geq\!1$. Note that we must have\\ 
\mbox{}\hfill  
$s\!=\!\min\limits_{j\in\{0,\ldots,k_{i,\zeta}-1\}} 
\left\{j+\ord_p\frac{f^{(j)}_{i-1,\mu}(\zeta_{i-1})}{j!}\right\}$, 
\hfill\mbox{}\\  
since $f_{i,\zeta}\!\in\!\left(\Z/\!\left(p^{k_{i,\zeta}}\right)\right)[x]$ 
for $i\!\geq\!1$. 
So then, the coefficient of $x^\ell$ in $f_{i-1,\mu}(\zeta_{i-1}+px)$ must 
be divisible by $p^{s+1}$ for all $\ell\!\geq\!s+1$. 
In other words, the coefficient of $x^\ell$ in $f_{i,\zeta}(x)$ must 
be divisible by $p$ for all $\ell\!\geq\!s+1$, and thus 
$\deg \tf_{i,\zeta}\!\leq\!s$. That $s\!\leq\!k_{i-1,\mu}-1$ follows 
from the definition of $s(f,\zeta)$, and $k_{i-1,\mu}\!\leq\!k$ since 
$k_{0,0}\!:=\!k$ and (thanks to Definition \ref{dfn:basic}) 
$k_{i-1,\mu}\!>\!k_{i,\zeta}$.

To prove the final bound, note that Lemma \ref{lemma:mult} implies that 
each term $s(f_{i-1,\mu},\zeta_{i-1})$ in the sum is 
at most the multiplicity of the root $\zeta_{i-1}$ of $\tf_{i-1,\mu}$. 
Since the sum of the multiplicities of the degenerate roots of 
$\tf_{i-1,\mu}$ is no greater than $\deg \tf_{i-1,\mu}$, 
we are done. \qed  

\medskip 
\noindent 
{\bf Assertion (4):} 
By Assertion (3), the sum of the degrees of the 
$\tf_{1,\zeta_0}$ (as $(f_{1,\zeta_0},k_{1,\zeta_0})$ ranges over all depth 
$1$ node labels of $\cT_{p,k}(f)$) is no greater than $\deg \tf_{0,0}$, which 
is at most $d$. 

By applying Assertion (3) to all nodes of depth $i\!\geq\!2$, 
the sum of the degrees of the $\tf_{i,\zeta}$ (as 
$(f_{i,\zeta},k_{i,\zeta})$ ranges over all depth $i$ node labels of 
$\cT_{p,k}(f)$) is no greater than the sum of the degrees of the 
$\tf_{i-1,\mu}$ (as $(f_{i-1,\mu},k_{i-1,\mu})$ ranges over all depth $i-1$ 
node labels of $\cT_{p,k}(f)$). 

Since $\deg \tf_{0,0}\!\leq\!d$ we thus obtain that, for {\em every} depth 
$i$, the sum of the degrees of the $\tf_{i,\zeta}$ (as 
$(f_{i,\zeta},k_{i,\zeta})$ ranges over all depth $i$ node labels of 
$\cT_{p,k}(f)$) is no greater than $d$. So by the final part of Assertion (3), 
our tree $\cT_{p,k}(f)$ has no more than $\floor{d/2}$ nodes at any fixed  
depth $\geq\!1$. So by Assertion (1) we are done. \qed

\medskip 
We are at last ready to prove our main theorem.

\medskip 
\noindent 
{\bf Proof of Theorem \ref{thm:tree}:} 
Since we already proved at the end of the last section 
that Algorithm \ref{algor:main} is correct, 
it suffices to prove the stated complexity bound for 
Algorithm \ref{algor:main}. Proving that Algorithm \ref{algor:main} 
runs as fast as stated will follow easily from (a) the fast 
randomized Kedlaya-Umans factoring algorithm from \cite{kedlaya} and  
(b) applying Lemmata \ref{lemma:ref} and \ref{lemma:branch} to show that 
the number of necessary factorizations and $p$-adic valuation calculations is 
well-bounded. 

More precisely, the {\bf For} loops and recursive calls of 
Algorithm \ref{algor:main} can be interpreted as a 
depth-first search of $\cT_{p,k}(f)$, with $\cT_{p,k}(f)$ being built 
along the way. In particular, we begin at the root node by factoring  
$\tf_{0,0}\!=\!\tf$ in $(\Z/(p))[x]$ via \cite{kedlaya}, in order to 
find the degenerate roots of $\tf$. (Factoring in fact dominates 
the complexity of the gcd computation that gives us $n_p(f_{0,0})$,  
if we use a deterministic near linear-time gcd algorithm such as that of 
Knuth and Sch\"onhage (see, e.g., \cite[Ch.\ 3]{bcs}).) 
This factorization takes time $\left(d^{1.5}\log p\right)^{1+o(1)}
+\left(d\log^2 p\right)^{1+o(1)}$ and requires $O(d\log p)$ random bits. 

Now, in order to continue the recursion, we need to compute $p$-adic 
valuations of polynomial coefficients in order to find the 
$s(f_{0,0},\zeta_0)$ and determine the edges emanating from our root. 
Expanding each $f_{0,0}(\zeta_0+px)$ can clearly be done mod $p^k$, so 
each such expansion takes time no worse than 
$d^2(k\log p)^{1+o(1)}$ via Horner's method and fast finite ring 
arithmetic (see, e.g., \cite{bs,vzgbook}). 
Computing $s(f_{0,0},\zeta_0)$ then takes time no worse than 
$d(k\log p)^{1+o(1)}$ using, say, the standard binary method for evaluating 
powers of $p$. There are no more than $\floor{d/2}$ possible 
$\zeta_0$ (thanks to Lemma \ref{lemma:branch}), so the total 
work so far is 
\[d^3(k\log p)^{1+o(1)}+\left(d\log^2 p\right)^{1+o(1)}.\] 
(To simplify our bound, we are rolling multiplicative constants 
into the exponent, at the price of a negligible increase in the 
little-$o(\cdot)$ terms in the exponent.) Note that now, 
computing the expansion $f_{0,0}(\zeta_0+px)$ dominates the 
factorization of $\tilde{f}_{0,0}$.

The remaining work can then be bounded similarly, but with one 
small twist: By Assertion (4) of 
Lemma \ref{lemma:branch}, the number of nodes at depth $i$ of our 
tree is never more than $\floor{d/2}$ and, by Lemma \ref{lemma:ref},  
the sum of the degrees of the $\tf_{i,\zeta}$ at level $i$ is 
no greater than $d$.  

Now observe that (for $i\!\geq\!2$) the amount of work needed to compute the 
$s(f_{i-1,\mu},\zeta_{i-1})$ at level $i-1$ (which are used to 
define the polynomials at level $i$) is 
no greater than\linebreak 
$d\cdot d(k\log p)^{1+o(1)}$, and this will be dominated by 
the subsequent computations of the expansions of the $f_{i,\zeta}$. In 
particular, by the 
basic calculus inequality $r^t_1+\cdots+r^t_\ell\!\leq\!(r_1+\cdots+
r_\ell)^t$ (valid for any $r,t\!\geq\!1$), the total amount of work for 
the factorizations for each subsequent level of $\cT_{p,k}(f)$ will be   
\[d^{1.5}(\log p)^{1+o(1)}+\left(d\log^2 p\right)^{1+o(1)},\] 
with $O(d\log p)$ random bits needed. 
The expansions of the $f_{i,\zeta}$ at level $i$ will take time 
no greater than $d^3(k\log p)^{1+o(1)}$ to compute. 
So our total work at each subsequent level is then 
\[ d^3(k\log p)^{1+o(1)}+\left(d\log^2 p\right)^{1+o(1)}.\] 

So then, 
the total amount of work for our entire tree will be 
\[ kd^3(k\log p)^{1+o(1)}+k\left(d\log^2 p\right)^{1+o(1)}.\] 
and the number of random bits needed is $O(dk\log p)$.  

We are nearly done, but we must still ensure that our algorithm has the 
correct Las Vegas properties. In particular, while finite field factoring 
can be assumed to succeed with probability $\geq\!2/3$, we use 
multiple calls to finite field factoring, {\em each of which could fail}.  
The simplest solution is to simply run our finite field factoring 
algorithm sufficiently many times to reduce the {\em over-all} 
error probability. In particular, thanks to Lemma \ref{lemma:branch}, 
and the basic union bound for probabilities, 
it is enough to enforce a success probability of $O\!\left(\frac{1}
{dk}\right)$ for each application of finite field factoring. 
This implies that we should run the algorithm from \cite{kedlaya} 
$O(\log(dk))$ many times each time we need a factorization over 
$(\Z/(p))[x]$. So, multiplying our last total by $\log(dk)$, 
this yields a final complexity bound of 
\[ kd^3(k\log p)^{1+o(1)}+\left(dk\log^2 p\right)^{1+o(1)}\] 
(since computing the expansions of the $f_{i,\mu}(\zeta_{i-1}+x)$ 
dominates our complexity) 
and a total number of $O(d k \log(dk)\log p)$ random bits needed. 

To conclude, note that as our algorithm proceeds with depth first 
search, we need only keep track of collections of $f_{i,\zeta}$  
occuring along a root-to-leaf path in $\cT_{p,k}(f)$. A polynomial 
of degree $d$ with integer coefficients all of absolute value less than 
$p^k$ requires $O(dk\log p)$ bits to store, and 
Lemma \ref{lemma:branch} tells us that the depth of $\cT_{p,k}(f)$ 
is $O(k)$. So we never need more than $O(dk^2\log p)$ bits of memory. \qed 

\section*{Acknowledgements} 
We are grateful to Qi Cheng, Shuhong Gao, and Daqing Wan for 
many useful conversations. We also humbly thank Joachim von zur Gathen 
for his kind encouragement, and the referee for suggestions that 
clarified our exposition. 

\bibliographystyle{amsalpha}

\end{document}